\documentclass[12pt]{article}
\usepackage{amssymb,amsmath,amsthm,amsfonts,enumerate, bbm, mathdots}

\usepackage{latexsym}
\usepackage{amscd}
\usepackage{stmaryrd}
\usepackage{color}
\usepackage[all]{xypic}
\usepackage{epsfig}
\usepackage{graphics}
\usepackage{ifthen}
\usepackage{varioref}
\usepackage{rotating}
\usepackage{extarrows}
\usepackage{cite}
\usepackage{mathrsfs}

\numberwithin{equation}{section}

\theoremstyle{plain}
\newtheorem{thm}{Theorem}[section]

\newtheorem{lem}[thm]{Lemma}

\definecolor{darkgreen}{rgb}{0.0625,0.64,0.0625}
\usepackage[
pdfstartview=FitH,
CJKbookmarks=true,
bookmarksnumbered=true,
bookmarksopen=true,
colorlinks,
pdfborder=001,
linkcolor=blue,
anchorcolor=green,
citecolor=red
]{hyperref}

\renewcommand{\Re}{\operatorname{Re}}
\newfont{\scyr}{wncyr10 scaled 550}

\def\proof{\noindent {\bf Proof.\;}}

\allowdisplaybreaks

\begin{document}
	
	\title{The reciprocals of tails of the alternating Riemann zeta function}
	
	\date{\small ~ \qquad\qquad School of Mathematical Sciences, Tongji University \newline No. 1239 Siping Road,
		Shanghai 200092, China}
	
	\author{Zhonghua Li ~and ~Lu Yan}

	\maketitle
	
	\begin{abstract}
		In this paper, we give the integer parts of reciprocals of tails of the alternating Riemann zeta function at $s=1,2,3,4$ by using several new inequalities and elementary method.
	\end{abstract}
	
	
	
	
	\section{Introduction}\label{Sec:Intro}
	
	In this paper, we study the reciprocals of tails of the alternating Riemann zeta function. The Riemann zeta function is defined by
	$$\zeta(s) =\sum_{k=1}^\infty\frac{1}{k^s},$$
	where $s \in \mathbb{C}$. The above series converges absolutely for $\Re(s)>1$, and it can be analytically  continued to the whole complex plane except for a  simple pole $s=1$ of residue $1$. The alternating Riemann zeta function is defined by
	$$\zeta^*(s) =\sum_{k=1}^\infty\frac{(-1)^{k+1}}{k^s},$$
	which converges for $\Re(s)>0$. As
	$$\zeta(s)=\frac{1}{1-2^{1-s}}\zeta^{\ast}(s)$$
	for $\Re(s)>1$, the above gives an explicit analytic continuation of the Riemann zeta function to the half plane $\Re(s)>0$ by the alternating Riemann zeta function.

	There are various elegant properties of the Riemann zeta function. Recently, some authors started to study the tails of the Riemann zeta function and the tails of the alternating Riemann zeta function, which are defined respectively by
	$$\zeta_n(s) =\sum_{k=n}^\infty\frac{1}{k^s},$$
	where $s \in \mathbb{C}$ with $\Re(s)>1$, and
	$$\zeta_n^*(s) =\sum_{k=n}^\infty\frac{(-1)^{k+1}}{k^s},$$
	where $s \in \mathbb{C}$ with $\Re(s)>0$. For example, in 2016 Lin \cite{Lin} studied the integer parts of reciprocals of tails of the Riemann zeta function at the integer point $s\geq 2$ and proved that
	$$\left\lfloor\zeta_n(2)^{-1}\right\rfloor=n-1,$$
	and
	$$\left\lfloor\zeta_n(3)^{-1}\right\rfloor=2n(n-1),$$
	where $n$ is any positive integer and $\lfloor x \rfloor$ denotes the largest integer less than or equal to $x$. Soon afterwards, Lin and Li\cite{Lin-Li} considered the computational formula for the case $s=4$, and they obtained:
	\[
	\left\lfloor\zeta_n(4)^{-1}\right\rfloor=
	\begin{cases}
		24m^3-18m^2+\left\lfloor\frac{3(5m-1)}{2}\right\rfloor &\text{if $n=2m$;}\\
		&\\
		24m^3-54m^2+\left\lfloor\frac{3(58m-17)}{4}\right\rfloor &\text{if $n=2m-1$,}
	\end{cases}
	\]
	for any positive integer $n$. Along this line, Xu \cite{Xu} also proved two computation formulas for $s=4,5$, and Hwang and Song \cite{Hwang-Song1} obtained a complicated formula for the case $s=6$, which depends on the residue of $n$ modulo $48$. In 2018, Kim and Song \cite{Kim-Song} studied the integer parts of the inverses of tails of the alternating Riemann zeta function for $s=\frac{1}{2},\frac{1}{3},\frac{1}{4}$, and they obtained that for any positive integer $n$ and $s=\frac{1}{2},\frac{1}{3},\frac{1}{4}$,
	\begin{equation}\label{1.1}
		\left\lfloor\zeta_n^{\ast}(s)^{-1}\right\rfloor=\left\lfloor (-1)^{n+1}2\left( n-\frac{1}{2}\right)^s\right\rfloor.
	\end{equation}
	Later, Hwang and Song \cite{Hwang-Song2} proved that when $s=\frac{1}{p}$ for any integer with $p\geq5$ or $s=\frac{2}{p}$ for any odd integer with $p\geq5$, there exists an integer $N>0$ such that the formula \eqref{1.1} still holds for every integer $n\geq N$.
	
	We want to mention that some similar questions on Fibonacci numbers have already been considered. For example, in 2008 Ohtsuka and Nakamura \cite{Ohtsuka-Nakamura} studied the properties of infinite sums of reciprocal Fibonacci numbers, and they proved that
	\[
	\left\lfloor\left(\sum_{k=n}^{\infty}\frac{1}{F_k}\right)^{-1}\right\rfloor=
	\begin{cases}
		F_{n-2} &\text{if $n\geq2$ is even;}\\
		F_{n-2}-1 &\text{if $n\geq1$ is odd,}
	\end{cases}
	\]
	and
	\[
	\left\lfloor\left(\sum_{k=n}^{\infty}\frac{1}{F_k^2}\right)^{-1}\right\rfloor=
	\begin{cases}
		F_{n-1}F_n-1 &\text{if $n\geq2$ is even;}\\
		F_{n-1}F_n &\text{if $n\geq1$ is odd.}
	\end{cases}
	\]
	Here, the Fibonacci sequence $\{F_k\}$ is defined by the recursive formula $F_{k+1}=F_{k}+F_{k-1}$ with $k\geq 1$ and the initial values $F_0=0$ and $F_1=1$. And then in 2013, Kuhapatanakul\cite{Kuhapatanakul} considered the infinite sums of reciprocal generalized Fibonacci numbers defined by the recursive formula $V_{k+1}=aV_{k}+bV_{k-1}$ and the initial values $V_0=c$ and $V_1=1$. And it was shown in \cite{Kuhapatanakul} that
	$$\left\lfloor\left(\sum_{k=n}^{\infty}\frac{(-1)^k}{V_k}\right)^{-1}\right\rfloor=(-1)^n(V_n+V_{n-1})-1,$$
	for any positive integer $n$, integers $a, b$ with $1\leq b \leq a$ and $c=0$. A direct corollary of the above formula is the following result for the Fibonacci numbers: 
	$$\left\lfloor\left(\sum_{k=n}^{\infty}\frac{(-1)^k}{F_k}\right)^{-1}\right\rfloor=(-1)^nF_{n+1}-1,$$
	where $n$ is any positive integer.
	
	Inspired by their results, we try to find the closed formulas of the integer parts of reciprocals of tails of the alternating Riemann zeta function at integer points. By using elementary method and several new inequalities, we obtain four interesting computational formulas for $\left\lfloor\zeta_n^*(s) ^{-1}\right\rfloor$ with $s=1,2,3,4$. That is, we shall prove the following four theorems.
	
	\begin{thm}\label{Thm:s=1}
		For a positive integer $n$, we have
		\[
		\left\lfloor\zeta_n^*(1) ^{-1}\right\rfloor=
		\begin{cases}
			-2n &\text{if $n\geq2$ is even;}\\
			2n-1 &\text{if $n\geq1$ is odd.}
		\end{cases}
		\]
	\end{thm}
	\begin{thm}\label{Thm:s=2}
		For a positive integer $n$, we have
		\[
		\left\lfloor\zeta_n^*(2) ^{-1}\right\rfloor=
		\begin{cases}
			-(2n^2-2n+1)-1 &\text{if $n\geq2$ is even;}\\
			2n^2-2n+1 &\text{if $n\geq1$ is odd.}
		\end{cases}
		\]
	\end{thm}
	\begin{thm}\label{Thm:s=3}
		For a positive integer $n$, we have
		\[
		\left\lfloor\zeta_n^*(3) ^{-1}\right\rfloor=
		\begin{cases}
			-\left( 2n^3-3n^2+\frac{9}{2}n-\frac{5}{2}\right)-\frac{3}{2} &\text{if $n\geq22$ is even;}\\
			2n^3-3n^2+\frac{9}{2}n-\frac{5}{2} &\text{if $n\geq7$ is odd.}
		\end{cases}
		\]
	\end{thm}
	\begin{thm}\label{Thm:s=4}
		For a positive integer $n$, we have
		\[
		\left\lfloor\zeta_n^*(4) ^{-1}\right\rfloor=
		\begin{cases}
			-(2n^4-4n^3+8n^2-6n-8)-1 &\text{if $n\geq10$ is even;}\\
			2n^4-4n^3+8n^2-6n-8&\text{if $n\geq11$ is odd.}
		\end{cases}
		\]
	\end{thm}
	
	The structure of this paper is as follows. In Section 2, we construct several inequalities which are necessary to the proofs of our theorems. In Section 3, we firstly prove Theorem \ref{Thm:foralls} which is a unified idea for all integers $s\geq1$, and then the proofs of  Theorems \ref{Thm:s=1} - \ref{Thm:s=4} are given respectively.
	
	
	\section{Several inequalites}\label{Sec:SIR}
	
	For a fixed integer $s$ with $1\leq s \leq4$, we want to prove there exist functions $f_s(k), g_s(k)$ with $\lim_{k \to \infty}{f_s(k)}=\lim_{k \to \infty}{g_s(k)}=\infty$ and integers $k_{s,even}, k_{s,odd}\geq1$ such that
	$$\frac{1}{f_s(k)+1}-\frac{1}{f_s(k+1)+1}<-\frac{1}{(2k)^s}+\frac{1}{(2k+1)^s}<\frac{1}{f_s(k)}-\frac{1}{f_s(k+1)}$$
	holds for any integer $k\geq k_{s,even}$, and
	$$\frac{1}{g_s(k)+1}-\frac{1}{g_s(k+1)+1}<\frac{1}{(2k-1)^s}-\frac{1}{(2k)^s}<\frac{1}{g_s(k)}-\frac{1}{g_s(k+1)}$$
	holds for any integer $k\geq k_{s,odd}$.
	
	\begin{lem}\label{Lem:s=1}
		Let $f_1(k)=-4k$ and $g_1(k)=4k-3.$
		For any positive integer $k$, we have
		\begin{align}\label{evens=1}
			\frac{1}{f_1(k)+1}-\frac{1}{f_1(k+1)+1}<-\frac{1}{2k}+\frac{1}{2k+1}<\frac{1}{f_1(k)}-\frac{1}{f_1(k+1)}
		\end{align}
		and
		\begin{align}\label{odds=1}
			\frac{1}{g_1(k)+1}-\frac{1}{g_1(k+1)+1}<\frac{1}{2k-1}-\frac{1}{2k}<\frac{1}{g_1(k)}-\frac{1}{g_1(k+1)}.
		\end{align}
	\end{lem}
	\proof
	As for any positive integer $k$ it holds
	$$-\frac{1}{2k}+\frac{1}{2k+1}+\frac{1}{4k}-\frac{1}{4k+4}=\frac{-1}{4k(k+1)(2k+1)}<0,$$
	the right-hand side of \eqref{evens=1} is proved. For the left-hand side of  \eqref{evens=1}, we use
	$$-\frac{1}{2k}+\frac{1}{2k+1}+\frac{1}{4k-1}-\frac{1}{4k+3}=\frac{3}{2k(2k+1)(4k-1)(4k+3)}>0 $$
	which holds for any positive integer $k$.
	
	Similarly, let's prove the inequality \eqref{odds=1}. The right-hand side follows from
	$$\frac{1}{2k-1}-\frac{1}{2k}-\frac{1}{4k-3}+\frac{1}{4k+1}=\frac{-3}{2k(2k-1)(4k-3)(4k+1)}<0, $$
	and the left-hand side follows from
	$$\frac{1}{2k-1}-\frac{1}{2k}-\frac{1}{4k-2}+\frac{1}{4k+2}=\frac{1}{2k(2k-1)(2k+1)}>0. $$
	This completes the proof of this lemma.
	\qed
	
	\begin{lem}\label{Lem:s=2}
		Let $f_2(k)=-2(2k-\frac{1}{2})^2-\frac{3}{2}$ and $g_2(k)=2(2k-\frac{3}{2})^2+\frac{1}{2}.$
		For any positive integer $k$, we have
		\begin{align}\label{evens=2}
			\frac{1}{f_2(k)+1}-\frac{1}{f_2(k+1)+1}<-\frac{1}{(2k)^2}+\frac{1}{(2k+1)^2}<\frac{1}{f_2(k)}-\frac{1}{f_2(k+1)}
		\end{align}
		and
		\begin{align}\label{odds=2}
			\frac{1}{g_2(k)+1}-\frac{1}{g_2(k+1)+1}<\frac{1}{(2k-1)^2}-\frac{1}{(2k)^2}<\frac{1}{g_2(k)}-\frac{1}{g_2(k+1)}.
		\end{align}
	\end{lem}
	\proof
	Since
	$$-\frac{1}{(2k)^2}+\frac{1}{(2k+1)^2}=-\frac{4k+1}{16k^4+16k^3+4k^2}$$
	and
	$$\frac{1}{f_2(k)}-\frac{1}{f_2(k+1)}=-\frac{16k+4}{64k^4+64k^3+16k^2+12},$$
	the right-hand side of \eqref{evens=2} is equivalent to
	$$\left(16k+4\right)\left(16k^4+16k^3+4k^2\right)<\left(4k+1\right)\left(64k^4+64k^3+16k^2+12\right).$$
	The above inequality is just $12>0$, which holds for any positive integer $k$. Using
	$$\frac{1}{f_2(k)+1}-\frac{1}{f_2(k+1)+1}=-\frac{16k+4}{64k^4+64k^3-8k+5},
	$$
	we find the left-hand side of \eqref{evens=2} is equivalent to
	$$\left(16k+4\right)\left(16k^4+16k^3+4k^2\right)>\left(4k+1\right)\left(64k^4+64k^3-8k+5\right).$$
	Through proper simplification, the above inequality is just
	$$16k^2+8k-5>0,$$
	which holds for any positive integer $k$. So this completes  the proof of the inequality \eqref{evens=2}.
	
	Similarly, for  the inequality \eqref{odds=2}, we have
	$$\frac{1}{(2k-1)^2}-\frac{1}{(2k)^2}=\frac{4k-1}{16k^4-16k^3+4k^2}.$$
	Since
	$$\frac{1}{g_2(k)}-\frac{1}{g_2(k+1)}=\frac{16k-4}{64k^4-64k^3+8k+5},$$
	through proper simplification, the right-hand side of \eqref{odds=2} is just
	$$16k^2-8k-5>0,$$
	which holds for any positive integer $k$. For the left-hand side, we have
	$$\frac{1}{g_2(k)+1}-\frac{1}{g_2(k+1)+1}=\frac{16k-4}{64k^4-64k^3+16k^2+12}.$$
	Then the left-hand side of \eqref{odds=2} is just $12>0$, which holds for any positive integer $k$. So this completes  the proof of the inequality \eqref{odds=2}.
	\qed

	\begin{lem}\label{Lem:s=3even}
		Let $f_3(k)=-16k^3+12k^2-9k+1.$
		For any positive integer $k\geq11$, we have
		\begin{align}\label{evens=3}
			\frac{1}{f_3(k)+1}-\frac{1}{f_3(k+1)+1}<-\frac{1}{(2k)^3}+\frac{1}{(2k+1)^3}<\frac{1}{f_3(k)}-\frac{1}{f_3(k+1)}.
		\end{align}
	\end{lem}
	\proof
	We have
	$$-\frac{1}{(2k)^3}+\frac{1}{(2k+1)^3}=-\frac{12k^2+6k+1}{64k^6+96k^5+48k^4+8k^3}.$$
	Since
	$$\frac{1}{f_3(k)}-\frac{1}{f_3(k+1)}=-\frac{48k^2+24k+13}{256k^6+384k^5+240k^4+104k^3+117k^2+75k-12},$$
	the right-hand side of \eqref{evens=3} is equivalent to
	\begin{align*}
		&\left(12k^2+6k+1\right)\left(256k^6+384k^5+240k^4+104k^3+117k^2+75k-12\right)>\\
		&\left(48k^2+24k+13\right)\left(64k^6+96k^5+48k^4+8k^3\right).
	\end{align*}
	The above inequality is just
	$$288k^5+1452k^4+1602k^3+423k^2+3k-12>0,$$
	which holds for all integers $k\geq1$.
	Since
	$$\frac{1}{f_3(k)+1}-\frac{1}{f_3(k+1)+1}=-\frac{48k^2+24k+13}{256k^6+384k^5+240k^4+72k^3+93k^2+33k-22},$$
	the left-hand side of \eqref{evens=3} is equivalent to
	\begin{align*}
		&\left(12k^2+6k+1\right)\left(256k^6+384k^5+240k^4+72k^3+93k^2+33k-22\right)<\\
		&\left(48k^2+24k+13\right)\left(64k^6+96k^5+48k^4+8k^3\right).
	\end{align*}
	The above inequality is just
	$$96k^5-972k^4-922k^3-27k^2+99k+22>0,$$
	which holds for all integers $k\geq11$. So this completes the proof.
	\qed

	\begin{lem}\label{Lem:s=3odd}
		Let $g_3(k)=16k^3-36k^2+33k-12.$
		For any positive integer $k\geq4$, we have
		\begin{align}\label{odds=3}
			\frac{1}{g_3(k)+1}-\frac{1}{g_3(k+1)+1}<\frac{1}{(2k-1)^3}-\frac{1}{(2k)^3}<\frac{1}{g_3(k)}-\frac{1}{g_3(k+1)}.
		\end{align}
	\end{lem}
	\proof
	We have
	$$\frac{1}{(2k-1)^3}-\frac{1}{(2k)^3}=\frac{12k^2-6k+1}{64k^6-96k^5+48k^4-8k^3}.$$
	Since
	$$\frac{1}{g_3(k)}-\frac{1}{g_3(k+1)}=\frac{48k^2-24k+13}{256k^6-384k^5+240k^4-104k^3+117k^2-75k-12},$$
	the right-hand side of \eqref{odds=3} is equivalent to
	\begin{align*}
		&\left(12k^2-6k+1\right)\left(256k^6-384k^5+240k^4-104k^3+117k^2-75k-12\right)<\\
		&\left(48k^2-24k+13\right)\left(64k^6-96k^5+48k^4-8k^3\right).
	\end{align*}
	The above inequality is just
	$$288k^5-1452k^4+1602k^3-423k^2+3k+12>0,$$
	which holds for all integers $k\geq4$.
	Since
	$$\frac{1}{g_3(k)+1}-\frac{1}{g_3(k+1)+1}=\frac{48k^2-24k+13}{256k^6-384k^5+240k^4-72k^3+93k^2-33k-22},
	$$
	the left-hand side of \eqref{odds=3} is equivalent to
	\begin{align*}
		&\left(12k^2-6k+1\right)\left(256k^6-384k^5+240k^4-72k^3+93k^2-33k-22\right)>\\
		&\left(48k^2-24k+13\right)\left(64k^6-96k^5+48k^4-8k^3\right).
	\end{align*}
	The above inequality is just
	$$96k^5+972k^4-922k^3+27k^2+99k-22>0,$$
	which holds for all integers $k\geq1$. This proves the inequality \eqref{odds=3}.
	\qed

	\begin{lem}\label{Lem:s=4even}
		Let $f_4(k)=-32k^4+32k^3-32k^2+12k+7.$
		For any positive integer $k\geq5$, we have
		\begin{align}\label{evens=4}
			\frac{1}{f_4(k)+1}-\frac{1}{f_4(k+1)+1}<-\frac{1}{(2k)^4}+\frac{1}{(2k+1)^4}<\frac{1}{f_4(k)}-\frac{1}{f_4(k+1)}.
		\end{align}
	\end{lem}
	\proof
	We have
	$$-\frac{1}{(2k)^4}+\frac{1}{(2k+1)^4}=-\frac{32k^3+24k^2+8k+1}{256k^8+512k^7+384k^6+128k^5+16k^4}.$$
	Since
	\begin{align*}
		&\frac{1}{f_4(k)}-\frac{1}{f_4(k+1)}=\\
		&-\frac{128k^3+96k^2+96k+20}{1024k^8+2048k^7+2048k^6+1280k^5+448k^4+64k^3-1488k^2-744k-91},
	\end{align*}
	the right-hand side of \eqref{evens=4} is equivalent to
	$$2048k^7+3584k^6-45312k^5-58880k^4-32608k^3-9624k^2-1472k-91>0,$$
	which holds for all integers $k\geq5$. Since
	\begin{align*}
		&\frac{1}{f_4(k)+1}-\frac{1}{f_4(k+1)+1}=\\
		&-\frac{128k^3+96k^2+96k+20}{1024k^8+2048k^7+2048k^6+1280k^5+384k^4-1648k^2-816k-96},
	\end{align*}
	the left-hand side of \eqref{evens=4} is equivalent to
	$$52480k^5+65600k^4+35840k^3+10480k^2+1584k+96>0,$$
	which holds for all integers $k\geq1$. This completes the proof.
	\qed

	\begin{lem}\label{Lem:s=4odd}
		Let $g_4(k)=32k^4-96k^3+128k^2-84k+12$.
		For any positive integer $k\geq6$, we have
		\begin{align}\label{odds=4}
			\frac{1}{g_4(k)+1}-\frac{1}{g_4(k+1)+1}<\frac{1}{(2k-1)^4}-\frac{1}{(2k)^4}<\frac{1}{g_4(k)}-\frac{1}{g_4(k+1)}.
		\end{align}
	\end{lem}
	\proof
	We have
	$$\frac{1}{(2k-1)^4}-\frac{1}{(2k)^4}=\frac{32k^3-24k^2+8k-1}{256k^8-512k^7+384k^6-128k^5+16k^4}.$$
	Since
	\begin{align*}
		&\frac{1}{g_4(k)}-\frac{1}{g_4(k+1)}=\\
		&\frac{128k^3-96k^2+96k-20}{1024k^8-2048k^7+2048k^6-1280k^5+384k^4-1648k^2+816k-96},
	\end{align*}
	the denominator of the above formula is greater than zero if $k\geq2$. Then for any integer $k\geq2$, the right-hand side of \eqref{odds=4} is equivalent to
	$$52480k^5-65600k^4+35840k^3-10480k^2+1584k-96>0,$$
	which holds for all integers $k\geq2$. Since
	\begin{align*}
		&\frac{1}{g_4(k)+1}-\frac{1}{g_4(k+1)+1}=\\
		&\frac{128k^3-96k^2+96k-20}{1024k^8-2048k^7+2048k^6-1280k^5+448k^4-64k^3-1488k^2+744k-91},
	\end{align*}
	the denominator of the above formula is greater than zero if $k\geq2$. Then for any integer $k\geq2$, the left-hand side of \eqref{odds=4} is equivalent to
	$$2048k^7-3584k^6-45312k^5+58880k^4-32608k^3+9624k^2-1472k+91>0,$$
	which holds for all integers $k\geq6$. This completes the proof.
	\qed


	\section{Proofs of the Theorems}\label{Sec:proof}
	For proving Theorems \ref{Thm:s=1}-\ref{Thm:s=4}, we give a unified idea of their proofs.
	\begin{thm}\label{Thm:foralls}
		Assume that for any positive integer $s\geq1$, there exist functions $f_s(x), g_s(x)\in \mathbb{Q}\left[ x\right] $, such that
		\begin{description}
			\item[(i)] $\lim_{k \to \infty}{f_s(k)}=\lim_{k \to \infty}{g_s(k)}=\infty$;
			\item[(ii)]  there exists a positive integer $k_{s,even}$, such that
			\begin{align}\label{evens}
				\frac{1}{f_s(k)+1}-\frac{1}{f_s(k+1)+1}<-\frac{1}{(2k)^s}+\frac{1}{(2k+1)^s}<\frac{1}{f_s(k)}-\frac{1}{f_s(k+1)}
			\end{align}
			holds for any integer $k\geq k_{s,even}$;
			\item[(iii)] there exists a positive integer $k_{s,odd}$, such that
			\begin{align}\label{odds}
				\frac{1}{g_s(k)+1}-\frac{1}{g_s(k+1)+1}<\frac{1}{(2k-1)^s}-\frac{1}{(2k)^s}<\frac{1}{g_s(k)}-\frac{1}{g_s(k+1)}
			\end{align}
			holds for any integer $k\geq k_{s,even}$.
		\end{description}	
		Then
		$$\frac{1}{f_s\left( \frac{n}{2}\right) +1}<\zeta_n^*(s)<\frac{1}{f_s\left( \frac{n}{2}\right) }$$
		holds for any positive even number $n$ with $n\geq 2k_{s,even}$, and
		$$\frac{1}{g_s\left( \frac{n+1}{2}\right) +1}<\zeta_n^*(s)<\frac{1}{g_s\left( \frac{n+1}{2}\right) }$$
		holds for any positive odd number $n$ with $n\geq 2k_{s,odd}-1$.
	\end{thm}
	
	\proof
	Let $n$ be a positive even integer. We have
	$$\zeta_n^*(s)=\sum_{k=\frac{n}{2}}^{\infty}\left(-\frac{1}{(2k)^s}+\frac{1}{(2k+1)^s}\right) .$$
	Using assumption (ii), we get
	\begin{align*}
		\sum_{k=\frac{n}{2}}^{\infty}\left( \frac{1}{f_s(k)+1}-\frac{1}{f_s(k+1)+1}\right)
		<&\sum_{k=\frac{n}{2}}^{\infty}\left(-\frac{1}{(2k)^s}+\frac{1}{(2k+1)^s}\right)\\
		<&\sum_{k=\frac{n}{2}}^{\infty}\left( \frac{1}{f_s(k)}-\frac{1}{f_s(k+1)}\right),
	\end{align*}
	which holds for any positive even number $n$ with $n\geq2k_{s,even}$. From the formulas above and assumption (i), we obtain
	$$\frac{1}{f_s\left( \frac{n}{2}\right) +1}<\zeta_n^*(s)<\frac{1}{f_s\left( \frac{n}{2}\right) }, $$
	which holds for any positive even number $n$ with $n\geq 2k_{s,even}$.
	
	Similarly, let $n$ be a positive odd number. We have
	$$\zeta_n^*(s)=\sum_{k=\frac{n+1}{2}}^{\infty}\left(\frac{1}{(2k-1)^s}-\frac{1}{(2k)^s}\right).$$
	Using assumption (iii), we obtain
	\begin{align*}
		\sum_{k=\frac{n+1}{2}}^{\infty}\left( \frac{1}{g_s(k)+1}-\frac{1}{g_s(k+1)+1}\right)
		<&\sum_{k=\frac{n+1}{2}}^{\infty}\left(\frac{1}{(2k-1)^s}-\frac{1}{(2k)^s}\right)\\
		<&\sum_{k=\frac{n+1}{2}}^{\infty}\left( \frac{1}{g_s(k)}-\frac{1}{g_s(k+1)}\right),
	\end{align*}
	which holds for any positive odd number $n$ with $n\geq2k_{s,odd}-1$. From the formulas above and assumption (i), we find
	$$\frac{1}{g_s\left( \frac{n+1}{2}\right) +1}<\zeta_n^*(s)<\frac{1}{g_s\left( \frac{n+1}{2}\right) },$$
	which holds for any positive even number $n$ with $n\geq 2k_{s,odd}-1$.
	\qed

	\subsection{Proof of Theorem \ref{Thm:s=1}}
	For $s=1$, we take $k_{1,even}=k_{1,odd}=1$, $f_1(k)=-4k$ and $g_1(k)=4k-3$. Then $f_1(k)$ and $g_1(k)$ satisfy the three conditions of Theorem \ref{Thm:foralls} by Lemma \ref{Lem:s=1}. Hence, using Theorem \ref{Thm:foralls}, we obtain that
	$$-\frac{1}{2n-1}<\zeta_n^*(1)<-\frac{1}{2n}$$
	holds for any positive even number $n$, and
	$$\frac{1}{2n}<\zeta_n^*(1)<\frac{1}{2n-1}$$
	holds for any positive odd number $n$. This completes the proof of Theorem \ref{Thm:s=1}.
	\qed

	\subsection{Proof of Theorem \ref{Thm:s=2}}
	For $s=2$, we take $k_{2,even}=k_{2,odd}=1$, $f_2(k)=-2(2k-\frac{1}{2})^2-\frac{3}{2}$ and $g_2(k)=2(2k-\frac{3}{2})^2+\frac{1}{2}$. Then using Lemma \ref{Lem:s=2} and Theorem \ref{Thm:foralls}, we obtain
	$$\frac{1}{f_2\left(\frac{n}{2}\right)+1}<\zeta_n^*(2)<\frac{1}{f_2\left(\frac{n}{2}\right)}$$
	or equivalently
	$$-\frac{1}{2n^2-2n+1}<\zeta_n^*(2)<-\frac{1}{2n^2-2n+2}$$
	holds for any positive even number $n$, and
	$$\frac{1}{g_2\left(\frac{n+1}{2}\right)+1}<\zeta_n^*(2)<\frac{1}{g_2\left(\frac{n+1}{2}\right)}$$
	or equivalently
	$$\frac{1}{2n^2-2n+2}<\zeta_n^*(2)<\frac{1}{2n^2-2n+1}$$
	holds for any positive odd number $n$. So this proves Theorem \ref{Thm:s=2}.
	\qed

	\subsection{Proof of Theorem \ref{Thm:s=3}}
	For $s=3$, we take $k_{3,even}=11$, $k_{3,odd}=4$, $f_3(k)=-16k^3+12k^2-9k+1$ and $g_3(k)=16k^3-36k^2+33k-12$. Then from Lemma \ref{Lem:s=3even}, Lemma \ref{Lem:s=3odd} and Theorem \ref{Thm:foralls}, we get
	$$\frac{1}{f_3\left(\frac{n}{2}\right)+1}<\zeta_n^*(3)<\frac{1}{f_3\left(\frac{n}{2}\right)},$$
	or equivalently
	$$-\frac{1}{	2n^3-3n^2+\frac{9}{2}n-2}<\zeta_n^*(3)<-\frac{1}{2n^3-3n^2+\frac{9}{2}n-1}$$
	holds for any positive even number $n\geq22$. So we have the inequality
	$$-\left( 2n^3-3n^2+\frac{9}{2}n-1\right) <\zeta_n^*(3) ^{-1}<-\left( 2n^3-3n^2+\frac{9}{2}n-2\right).$$
	Since $	-\left( 2n^3-3n^2+\frac{9}{2}n-1\right) $ and $-\left( 2n^3-3n^2+\frac{9}{2}n-2\right)$ are two consecutive integers, it follows that for any positive even number $n\geq22$,
	$$\left\lfloor\zeta_n^*(3) ^{-1}\right\rfloor=	-\left( 2n^3-3n^2+\frac{9}{2}n-1\right).$$
	Similarly, for any positive odd number $n\geq7$, we get
	$$\frac{1}{g_3\left(\frac{n+1}{2}\right)+1}<\zeta_n^*(3)<\frac{1}{g_3\left(\frac{n+1}{2}\right)},$$
	or equivalently
	$$\frac{1}{2n^3-3n^2+\frac{9}{2}n-\frac{3}{2}}<\zeta_n^*(3)<\frac{1}{2n^3-3n^2+\frac{9}{2}n-\frac{5}{2}}.$$
	Then we have the inequality
	$$2n^3-3n^2+\frac{9}{2}n-\frac{5}{2}<\zeta_n^*(3)^{-1}<2n^3-3n^2+\frac{9}{2}n-\frac{3}{2}.$$
	Since $2n^3-3n^2+\frac{9}{2}n-\frac{3}{2}$ and $2n^3-3n^2+\frac{9}{2}n-\frac{5}{2} $ are two consecutive positive integers, it follows that for any positive odd number $n\geq7$,
	$$\left\lfloor\zeta_n^*(3)^{-1}\right\rfloor=2n^3-3n^2+\frac{9}{2}n-\frac{5}{2}.$$
	
	This completes the proof of Theorem \ref{Thm:s=3}.
	\qed

	\subsection{Proof of Theorem \ref{Thm:s=4}}
	
	For $s=4$, we take $k_{4,even}=5$, $k_{4,odd}=6$, $f_4(k)=-32k^4+32k^3-32k^2+12k+7$ and $g_4(k)=32k^4-96k^3+128k^2-84k+12$. Hence by Lemma \ref{Lem:s=4even}, Lemma \ref{Lem:s=4odd} and Theorem \ref{Thm:foralls}, we have
	$$\frac{1}{f_4\left(\frac{n}{2}\right)+1}<\zeta_n^*(4)<\frac{1}{f_4\left(\frac{n}{2}\right)},$$
	or equivalently
	$$-\frac{1}{2n^4-4n^3+8n^2-6n-8}<\zeta_n^*(4)<-\frac{1}{2n^4-4n^3+8n^2-6n-7}$$
	holds for any positive even number $n\geq10$. Therefore, we find
	$$\left\lfloor\zeta_n^*(4)^{-1}\right\rfloor=-\left( 2n^4-4n^3+8n^2-6n-7\right) $$
	holds for any positive even number $n\geq10$.
	Similarly, for any positive odd number $n\geq11$, we have
	$$\frac{1}{g_4\left(\frac{n+1}{2}\right)+1}<\zeta_n^*(4)<\frac{1}{g_4\left(\frac{n+1}{2}\right)},$$
	or equivalently
	$$\frac{1}{2n^4-4n^3+8n^2-6n-7}<\zeta_n^*(4)<\frac{1}{2n^4-4n^3+8n^2-6n-8}.$$
	Then it follows that for any positive odd number $n\geq11$,
	$$\left\lfloor\zeta_n^*(4)^{-1}\right\rfloor=2n^4-4n^3+8n^2-6n-8.$$
	
	This completes the proof of Theorem \ref{Thm:s=4}.
	\qed


\end{document}